\documentclass[11pt]{amsart}
\usepackage{amssymb,latexsym,comment,url,amsmath}

\usepackage[T1]{fontenc}
\usepackage{rotating,graphicx}
\usepackage{currvita}
\usepackage{hyperref}
\usepackage{diagbox}

\usepackage{longtable}

\newcommand{\Char}{\operatorname{char}}
\newcommand{\Jac}{\operatorname{Jac}}

\newcommand{\Gal}{\operatorname{Gal}}
\newcommand{\Div}{\operatorname{Div}}
\newcommand{\GL}{\operatorname{GL}}
\newcommand{\divv}{\operatorname{div}}

\newcommand{\isom}{ \cong }

\newcommand{\PP}{{\mathbb P}}
\newcommand{\Q}{{\mathbb Q}}

\newcommand{\Z}{{\mathbb Z}}

\newenvironment{Proof}{\par\noindent{\sc Proof:}}%
                      {\hspace*{\fill}\nobreak$\Box$\par\medskip}
                       {\hspace*{\fill}\nobreak$\Box$\par\medskip}

\newtheorem{Proposition}{Proposition}[section]
\newtheorem{Theorem}[Proposition]{Theorem}

\newtheorem{Corollary}[Proposition]{Corollary}

\theoremstyle{definition}
\newtheorem{Definition}[Proposition]{Definition}
\newtheorem{Remark}[Proposition]{Remark}
 \newtheorem{Example}[Proposition]{Example}

\renewcommand{\baselinestretch}{1.1}

\begin{document}

\title[Quartic models and $D(q)$-quintuples]{Divisibility by 2 on quartic models of elliptic curves and rational Diophantine $D(q)$-quintuples}

\author[M. Sadek]%
{Mohammad~Sadek}
\address{Faculty of Engineering and Natural Sciences, Sabanc{\i} University, Tuzla, \.{I}stanbul, 34956 Turkey}
\email{mmsadek@sabanciuniv.edu}

\author[T. Yes\.{I}n]%
{Tu\u{g}ba Yes\.{I}n}

\email{tugbayesin@sabanciuniv.edu}

\date{}

\let\thefootnote\relax

\begin{abstract}
	Let $C$ be a smooth genus one curve described by a quartic polynomial equation over the rational field $\Q$ with $P\in C(\Q)$.  We give an explicit criterion for the divisibility-by-$2$ of a rational point on the elliptic curve $(C,P)$.  This provides an analogue to the classical criterion of the divisibility-by-$2$ on elliptic curves described by Weierstrass equations.  
	
We employ this criterion to investigate the question of extending a rational $D(q)$-quadruple to a quintuple.  We give concrete examples
to which we can give an affirmative answer.  One of these results implies that although the rational $ D(16t+9) $-quadruple $\{t, 16t+8,2 25t+14, 36t+20 \}$ can not be extended to a polynomial $ D(16t+9) $-quintuple using a linear polynomial,  there are infinitely many rational values of $t$ for which the aforementioned rational $ D(16t+9) $-quadruple can be extended to a rational $ D(16t+9) $-quintuple.  Moreover, these infinitely many values of $t$ are parametrized by the rational points on a certain elliptic curve of positive Mordell-Weil rank. 
\end{abstract}

\maketitle

\section{Introduction}

Let $ E $ be an elliptic curve over a number field $ K$.  The Mordell-Weil Theorem asserts that the abelian group of rational points $E(K)$ is finitely generated. In particular,  there are finitely many points $P_1,...,P_n$ in $E(K)$ such that any $P\in E(K)$ can be written as a linear combination $m_1P_1+\cdots+m_nP_n$ for some integers $m_1,\cdots,m_n$. During the course of the proof of the latter theorem, one proves the Weak Mordell-Weil Theorem which states that the abelian group $E(K)/2E(K)$ is finite. 

In order to show that $E(K)/2E(K)$ is finite, one needs to pass from a point $P\in E(K)$ to a point $Q\in E$ such that $2Q=P$. This process is the so-called $2$-descent on elliptic curves. The main ingredient to establish the $2$-descent is to give a necessary and sufficient condition such that $Q\in E(K)$.  In fact, one knows that if $E$ is described by the following Weierstrass equation 
\begin{equation*}
E: y^2=(x-a_{1})(x-a_{2})(x-a_{3}),\qquad a_i\textrm{ are distinct in }K,
\end{equation*}
then a point $ P=(x_{0}, y_{0})\in E(K) $ is divisible by 2 in $ E(K) $ if and only if all three elements $ x_{0}-a_{i} $ are squares in $ K $, see \cite[Chapter IV]{Knapp},  \cite[Chapter 6]{Husemoller}, \cite[Chapter VIII]{Silverman1}, or \cite{Bekker} for a criterion of the divisibility of rational points by powers of $2$. 

 The following quartic equation
\begin{equation*}
y^2= (a_1x+b_1)(a_2x+b_2)(a_3x+b_3)(a_4x+b_4), \qquad a_i,b_i\in K,
\end{equation*}
where $ b_{i}/a_{i} $ are distinct in $K$, describes a genus one curve $C$. Fixing a rational point $P\in C(K) $ to serve as the identity element of the group law, one may look for a similar criterion for the divisibility-by-$2$ on the elliptic curve $(C,P)$. It turns out that one can obtain a similar condition that depends on $P$, more precisely, a point $Q\in C(\Q)$ is twice a rational point if and only if the values of certain degree-$2$ polynomials evaluated at the $x$-coordinate of $Q$ are all squares in $\Q$.  Consequently, we show how to characterize such quartic models that possess rational $4$-torsion points.

A rational $D(q)$-$m$-tuple is an $m$-tuple $a_1, \cdots, a_m$ of distinct nonzero rational numbers such that
$a_ia_j + q$ is a square for all $1\le  i < j \le  m$.  If $q=1$, then the latter $m$-tuple is called a rational Diophantine $m$-tuple.The divisibility-by-$2$ on elliptic curves described by Weierstrass equations has been used to study rational $D(q)$-$m$-tuples. In \cite{Dujella-3}, $2$-divisibility on elliptic curves described by Weierstrass equations was used to extend  rational Diophantine triples to quadruples. It was also used to show that there are infinitely many rational Diophantine sextuples, see \cite{Dujella-6}. 
In \cite{Drazic},  it was proved that assuming the Parity Conjecture for the twists of several explicitly given elliptic curves, the density of rational numbers $q$
 for which there exist infinitely many rational $D(q)$-quintuples is at least $295026/296010 \approx 99.5\%$.

Rational Diophantine tuples have turned out to provide a useful tool to construct elliptic curves with prescribed torsion subgroups and high rank. In \cite{Dujella2021a}, rational Diophantine triples have been used to construct elliptic curves over $\Q(u)$ with rank 2 and either torsion subgroup $\Z/8\Z$ or $\Z/2\Z\times\Z/6\Z$. In \cite{DujellaSoydan}, for each of the
groups $\Z/2\Z \times \Z/k\Z$ for $k = 2, 4, 6, 8$, the authors proved the existence of infinitely many rational Diophantine quadruples with the property that the induced elliptic curve has this
torsion group. In \cite{Dujella2000}, the so-called regular Diophantine quadruples and quintuples, were characterized by elliptic curves. In addition, these characterizations were used to find examples
of elliptic curves over $\Q$ with torsion group $\Z/2\Z \times \Z/2\Z$ and with Mordell-Weil
rank equal $8$.

Researchers have been investigating $D(q)$-tuples whose elements enjoy certain properties. For example, in \cite{Dujella2021} the authors prove the existence of infinitely many essentially different $D(q)$-quintuples, where $q$ is an integer, whose elements are squares.
Further, integers that possess the $D(q)$-property for at least two integers $q_1$, $q_2$ have been studied. In fact, the authors of \cite{Dujella2020} proved the existence of infintely many essentially different sets consisting of perfect squares which are simultaneously $D(q_1)$-quadruples and $D(q_2)$-quadruples for distinct nonzero perfect squares $q_1$ and $q_2$.

For every rational number $q$, the authors of \cite{Goran} found all rational $m$ such
that there exists a rational $D(q)$-quadruple $\{a_1,a_2,a_3,a_4\}$ with product $a_1a_2a_3a_4 = m$.
 Using a certain rational map defined on a specific elliptic curve, the authors show that all such quadruples are identified if a certain rational map defined on the elliptic curve attains rational square values. For this reason, using the divisibility by-$2$ criterion on elliptic curves described by quartic equations in our work resembles the approach used in the aforementioned paper.

The notion of a strong rational $m$-tuple was introduced in \cite{Dujella2008}. Such a tuple is a rational Diophantine $m$-tuple, $\{a_1,\ldots,a_m\}$, with the additional property that $a_i^2+1$ is a rational
square for every $i = 1, \ldots , m$. The authors proved that there exist infinitely many strong rational Diophantine triples.  A strong rational $D(q)$-$m$-tuple
is a set of non-zero rationals $\{a_1, \ldots, a_m\}$ such that $a_ia_j +q$ is a square for all $i, j = 1, \ldots , m$, including the case $i = j$. The case $q = -1$ was studied
in \cite{Dujella2018} and it was shown that there exist infinitely many strong rational $D(-1)$-triples. In \cite{Sadek}, it was proved that there exist infinitely many square-free integers $q$ with the property that there
exist infinitely many strong rational $D(q)$-triples.

  A natural question to pose is how large a set of rational numbers enjoying the $D(q)$-property, for some $q\in \Q$, can be.  For a historical overview of rational $D(q)$-$m$-tuples, we refer the reader to \cite{Historical-1}, \cite[Sections 14.6 and 16.7]{Historical-2} as well as the webpage of Andrej Dujella\textsuperscript{1}\footnote{\textsuperscript{1}https://web.math.pmf.unizg.hr/$\sim$duje/dtuples.html}.

Jones initiated the study of polynomial $ D(q) $-$m$-tuples where $ q $ itself is a polynomial, see \cite{Jones1,Jones2}. If we define
\begin{equation*}
P_{q}=\sup \{ |S| : S\text{ is a polynomial $D(q)$-tuple} \},
\end{equation*}
then $ P_{q}\leq 22 $ for all $ q\in \mathbb{Z} $, see \cite[Theorem 1]{Dujella-7}. More properties of $ P_{q} $ can be found in \cite{Pq}. In this article, we focus on the case when $q$ is a linear polynomial. Setting
\begin{equation*}
L=\sup\{|S|: S \textrm{ is a polynomial $D(ax+b)$-tuple for some }a\ne 0 \textrm{ and }b\},
\end{equation*}
one can easily observe that $L\ge 4$ by viewing the $D(16x+9)$-quadruple $T=\{x, 16x+8, 25x+14, 36x+20 \} $, see \cite{polynomialformulas}. For upper bounds on $ L $, the reader can consult the papers \cite{Dujella,bounds}. 

We examine the case of polynomial $D(q)$-$m$-tuples consisting of linear polynomials where $q$ itself is a linear polynomial. If the set $ S $ consists only of linear polynomials, then $ \sup\{|S|\} $ is 4, see \cite{Dujella}. Hence, the above $D(16x+9)$-quadruple $ T  $ can not be extended to a polynomial $ D(16x+9) $-quintuple using a linear polynomial. In this work, we use the $2$-divisibility criterion on elliptic curves described by quartic models to study the extension of polynomial $D(q)$-qudruples to rational quintuples at infinitely many values of these polynomials. 
In fact, we show that although $T$ cannot be extended to a polynomial $D(16x+9)$-quintuple, there are infinitely many values for $x$ parametrized by an elliptic curve of positive rank such that $T$ can be extended to a quintuple using a rational function. We also present other polynomial $D(ax+b)$-quadruples with the same property. 

We remark that the theory of elliptic curves was used to show that there are only finitely many ways of extending a rational
$D(q)$-quadruple to a rational $D(q)$-quintuple, see \cite{Herrmann}. Our method provides an explicit description of how to extend certain rational $D(q)$-quadruples to rational $D(q)$-quintuples. In \cite{ondiophantine}, an explicit expression for the element extending a rational $D(q)$-quadruple to a raional $D(q)$-quintuple was provided if $q$ is a rational square. This means that if $x$ is chosen such that $16x+9$ is a rational square $q^2$, then our result together with \cite{ondiophantine} provides a method of constructing almost rational $D(q^2)$-sextuples, i.e.,  a tuple $a_1, \cdots, a_6$ of distinct nonzero rational numbers such that
$a_ia_j + q^2$ is a square for all $1\le  i < j \le  6$ except for $(i,j)=(5,6)$.  

\subsection*{Acknowledgments} The authors would like to thank the anonymous referees for several corrections, comments and valuable suggestions that improved the manuscript. They also would like to thank Andrej Dujella for several stimulating discussions. M. Sadek is partially supported by BAGEP Award of the Science Academy, Turkey.

\section{Models of elliptic curves}

In this section we introduce the genus one curve models that we are going to use throughout this paper.  
\subsection{Quartic models}
We recall that a {\em Weierstrass equation} is an equation of the form 
$$y^2+a_1xy+a_3y=x^3+a_2x^2+a_4x+a_6$$ where the coefficients $a_1,\cdots,a_6$ are lying in a field $K$.  One may associate to such equation the invariants $c_4$, $c_6$ and $\Delta$ which are polynomials in $a_1,\cdots, a_6$ with integer coefficients satisfying $1728\Delta=c_4^3-c_6^2$, \cite[Chapter III]{Silverman1}. If $\Delta\ne 0$, then the Weierstrass equation describes a smooth projective genus one curve with a $K$-rational point at infinity on the curve, i.e., an elliptic curve.  Two such equations describe the same curve if they are related via a transformation of the form $$x\mapsto u^2x+r,\qquad y\mapsto u^3y+u^2sx+t$$ for some $u\in K^{\times}$, $r,s,t\in K$. Therefore, if $\Char K\ne 2,3$, a Weierstrass equation can be written as $y^2=x^3+Ax+B$ for some $A,B\in K$.

\begin{Definition}
A {\em quartic model} is an equation of the form $y^2+P(x)y=Q(x)$ where $P$ and $Q$ are polynomials of degree $2$ and $4$, respectively, with coefficients in $K$. 

If $Q(x)+P(x)^2/4=ax^4+bx^3+cx^2+dx+e$, then we attach to the quartic model the invariants $c_4=2^4I$ and $c_6=2^5J$ where 
\begin{eqnarray*}
I&=&12ae-3bd+c^2,\\
J&=& 72ace-27ad^2-27b^2e+9bcd-2c^3.
\end{eqnarray*}
Moreover, the discriminant $\Delta$ of the quartic model satisfies $1728\Delta=c_4^3-c_6^2$.
\end{Definition}
If the discriminant $\Delta$ of a quartic model is different from $0$, then the quartic model describes a genus one curve $C$ over $K$. If the set of $K$-rational points $C(K)$ of $C$ is non-empty, then the quartic model describes an elliptic curve.
Two such models describe the same curve if they are related via a transformation of the form 
$$x\mapsto (a_{11}x+a_{21})/(a_{12}x+a_{22}),\qquad y\mapsto  \mu y+rx^2+sx+t$$ where $(a_{ij})\in \GL_2(K)$, $\mu\in K^{\times}$, $r,s,t\in K$. It follows that if $\Char K\ne 2,3$, then a quartic model can be written in the form $y^2=P(x)$ where $P(x)$ is a polynomial of degree $4$ with coefficients in $K$. 

If $C$ is a smooth genus one curve over a field $K$, $\Char K\ne 2,3$, defined by a quartic model, then the Jacobian $J$ of $C$ has Weierstrass equation $y^2=x^3-27c_4x-54c_6$, where $c_4$ and $c_6$ are the invariants of the quartic model, \cite[Theorem 2.8]{Fisher1}. In fact, if $C(K)\ne\emptyset$, then $C$ and $J$ are isomorphic over $K$.  

\subsection{Group law on quartic models}

Let $K$ be a field. The set of $K$-rational points on an elliptic curve $E$ defined over $K$ is an abelian group. If the curve is given by a Weierstrass equation, then the group law is described using the chord and tangent process. Given a smooth genus one curve $C$ that possesses a $K$-rational point and defined by a quartic model, one uses the isomorphism between the curve $C$ and its Jacobian to define the group law on $C$.

We consider a quartic model $y^2=f(x)$ where $f(x)\in K[x]$, $\deg f(x)=4$, and the discriminant of the model is nonzero, or equivalently, $f(x)$ has no multiple roots. This quartic model describes a smooth genus one curve over $K$. The Jacobian of $C$ will be denoted by $E$.  

From now on, we assume that $C(K)\ne\emptyset$.  This allows us to assume that the leading coefficient (or the constant term) of $f(x)$ is a square in $K$, \cite[Proposition 4.1]{Cremona1}. In this case, if the leading coefficient of $f(x)$ is a square $a^2$, $a\in K$, we set $\infty_+$ and $\infty_-$ to be the two rational points at infinity, namely, $(x:y:z)=(1:a:0)$ and $(1:-a:0)$, respectively.

 We fix a point $ P\in C(K) $. Let $\phi_P$ be a $K$-birational isomorphism between $C$ and $E$
 \begin{eqnarray*}
\phi_{P}: C \longrightarrow E \ \ \textrm{ such that } \ \phi_{P}(P)=O_{E}.
\end{eqnarray*}
	The map $\phi_P$ may be used to define an abelian group structure on $ C $ as follows 
	\begin{eqnarray*}
	Q_{1}+_{P}Q_{2}= \phi_{P}^{-1}(\phi_{P}(Q_{1})+\phi_{P}(Q_{2}))
	\end{eqnarray*}
	with $ P $ the identity on $(C,+_P)$.

 In particular, we say that $S\in n C(K)$, $n\ge 2$, if $ S=\underbrace{Q+_{P}\cdots +_P Q}_{n-\textrm{times}}$ for some $ Q\in C(K) $. This identifies $nC(K)$ with $nE(K)$.

\section{$2$-Divisibility on quartic models}
	
	Let $K$ be a perfect field of characteristic different from $2$. In this section, we consider quartic models $y^2=f(x)$ where $f(x)$ is a polynomial of degree $4$ with coefficients in $K$ and no multiple roots. We assume moreover that $f(x)$ splits completely in $K$. In other words, a quartic model will be of the form 
	$$y^2=f(x):=(a_1x+b_1)(a_2x+b_2)(a_3x+b_3)(a_4x+b_4),\qquad a_i\in K^{\times},b_i\in K,\; (a_ix+b_i)/(a_jx+b_j)\not\in K\textrm{ for }i\ne j,$$ and it describes a smooth genus one curve $C$ over $K$. The existence of the points $Q_i:=(-b_i/a_i,0)\in C(K)$, $1\le i\le 4$, imlies that $C(K)\ne\emptyset$.  We notice that $\infty_+,\infty_-\in C(K(\sqrt{a_1a_2a_3a_4}))$.
	We set $f_i(x):=a_ix+b_i$ and $c_i=-b_i/a_i$, $1\le i\le 4$.  
	
	We will always assume the existence of a rational point $(x_0,y_0)\in C(K)$ different from the points $Q_i$, $1\le i\le 4$. In particular, $|C(K)|>4 $. We fix throughout a $K$-birational isomorphism $\phi:C\to E:=\Jac(C) $ such that 
	 $\phi(x_0,y_0)=O_E$.
	 
	Now, we define the following rational maps $g_{ij}\in K(E)$, $ 1\le i, j\le 4$, as follows
	
	\[g_{ij}(P)=
      f_i(x_0)f_j(x_0)f_i(x(\phi^{-1}(P)))f_j(x(\phi^{-1}(P)))
\] 
	where $x(\phi^{-1}(P))$ is the $x$-coordinate of $\phi^{-1}(P)\in C(K)$.  
	It is clear that $g_{ij}(O_E)=f_i(x_0)^2f_j(x_0)^2\in (K^{\times})^2$.
	
	The following two propositions give properties of the maps $g_{ij}$ that we are going to use during the course of the proof of the main theorem of this section.  These properties have been proved for other rational maps on different models of elliptic curves, see for example \cite[Chapter 6]{Husemoller}, \cite[Chapter IV]{Knapp}, and \cite{Goran}. 
	\begin{Proposition}
	\label{Prop1} Let $[2]:E\to E$ be the multiplication by-$2$-morphism on $E$. 
There exist $h_{ij}\in K(E)$ such that $g_{ij}\circ [2]=h_{ij}^2$ for all $i,j$.
	\end{Proposition}
	\begin{Proof}
	One can see that $ \divv(g_{ij})=2\phi(Q_{i})+2\phi(Q_{j})-2\phi(\infty_{+})-2\phi(\infty_{-}) $.

 We set $[2]^*:\Div(E)\to \Div(E)$ to be the map $(Q)\mapsto \sum_{P\in[2]^{-1}Q}(Q)$.
Let $ \tilde{h}_{ij}\in \overline{K}(E) $ be such that
\begin{eqnarray*}
\divv(\tilde{h}_{ij})&=&[2]^*(\phi(Q_{i})+\phi(Q_{j})-\phi(\infty_{+})-\phi(\infty_{-}))\\
&=& \sum_{T\in E[2]}^{}(M_{i}+T)+ \sum_{T\in E[2]}^{}(M_{j}+T)- \sum_{T\in E[2]}^{}(N_{1}+T)- \sum_{T\in E[2]}^{}(N_{2}+T)
\end{eqnarray*}
where $ 2M_{i}=\phi(Q_{i}), \ 2M_{j}=\phi(Q_{j}), \ 2N_{1}=\phi(\infty_{+}), \ 2N_{2}=\phi(\infty_{-}) $. Then we can  observe that
\begin{eqnarray*}
 \divv(g_{ij}\circ [2])= 2\divv(\tilde{h}_{ij})=\divv(\tilde{h}_{ij}^{2}).
\end{eqnarray*}   
There exists $ r\in \overline{K} $ such that $ r\tilde{h}_{ij}^{2}= g_{ij}\circ [2]   $,  see for example \cite[Theorem 7.8.3]{Galbraith}. We define $ h_{ij}= \tilde{h}_{ij}\sqrt{r}  $.

 It is clear that $ \tilde{h}_{ij}\in K(E) $ for all $ i,j $. This follows by choosing $ \sigma \in \Gal(\overline{K}/K) $ and observing that $\sigma$ permutes the zeros of $ \tilde{h}_{ij} $, and the poles of $ \tilde{h}_{ij} $, respectively. More precisely,
\begin{eqnarray*}
O_{E}=  (\phi(Q_{i}))^{\sigma}-\phi(Q_{i})= (2M_{i})^{\sigma}-2M_{i}=2(M_{i}^{\sigma}-M_{i}),
\end{eqnarray*}
hence $ M_{i}^{\sigma}=M_{i}+T $ where $ T\in E[2] $. Same holds if one replaces $M_i$ with $N_i$. It is left to show that $r\in (K^{\times})^{2}$. This holds by evaluating both sides of the equality $ r\tilde{h}_{ij}^{2}= g_{ij}\circ [2]   $ at $O_E$. The statement holds as $g_{ij}(O_E)\in (K^{\times})^2$.
\end{Proof}
	\begin{Proposition}
	\label{Prop2}
	For any $P,Q\in E(K)$, one has
\[
g_{ij}(P+Q) \equiv g_{ij}(P)g_{ij}(Q)\mod K^2.
\]
	\end{Proposition}
\begin{Proof}
When $i=j$, the statement is straightforward, so we may pick $i\ne j$ and set $g:=g_{ij}$. According to Proposition \ref{Prop1}, we see that $g\circ [2]=h^2$ for some $h\in K(E)$. 

	Let $ P=2\tilde{P}, \ Q=2\tilde{Q} $. First we will prove 
	\begin{equation}\label{first prove}
	\dfrac{(h(\tilde{P}+\tilde{Q}))^{\sigma}}{h(\tilde{P}+\tilde{Q})}=\dfrac{(h(\tilde{P}))^{\sigma}}{h(\tilde{P})} \dfrac{(h(\tilde{Q}))^{\sigma}}{h(\tilde{Q})}
	\end{equation}
	for every $ \sigma \in \Gal(\bar{K}/K) $. Fix $T\in E[2] $, we have $ h^{2}(S+T)=g\circ[2](S+T)=g\circ[2](S)=h^{2}(S) $ for any $S\in E$, so $ \frac{h(S+T)}{h(S)}= \pm1 $. Considering the morphism $E\to \PP^1$ induced by the rational map $S\mapsto h(S+T)/h(S)$, one then may assume that it must be a constant map. Since $ 2\tilde{P}=P\in E(K),\ \  2\tilde{Q}=Q\in E(K) $, we get $ \tilde{P}^{\sigma}-\tilde{P}\in E[2], \,\tilde{Q}^{\sigma}-\tilde{Q}\in E[2], \ \  (\tilde{P}+\tilde{Q})^{\sigma}-(\tilde{P}+\tilde{Q})\in E[2] $ for every $ \sigma \in \Gal(\bar{K}/K) $. Now we have
	\begin{equation*}
	\dfrac{(h(\tilde{P}))^{\sigma}}{h(\tilde{P})}=\dfrac{h(\tilde{P}^{\sigma})}{h(\tilde{P})}=\dfrac{h(\tilde{P}+(\tilde{P}^{\sigma}-\tilde{P}))}{h(\tilde{P})}=\dfrac{h(S+(\tilde{P}^{\sigma}-\tilde{P}))}{h(S)}\quad \textrm{ for any } S\in E
	\end{equation*}
	Similarly,
	
	\begin{equation*}
	\dfrac{(h(\tilde{Q}))^{\sigma}}{h(\tilde{Q})}=\dfrac{h(S+(\tilde{Q}^{\sigma}-\tilde{Q}))}{h(S)}, \quad \textrm{ and }  \dfrac{(h(\tilde{P}+\tilde{Q}))^{\sigma}}{h(\tilde{P}+\tilde{Q})}=\dfrac{h(S+(\tilde{P}+\tilde{Q})^{\sigma}-(\tilde{P}+\tilde{Q}))}{h(S)}\qquad \textrm{ for any } S\in E.
	\end{equation*}
	Therefore,
	\begin{eqnarray*}
		\dfrac{(h(\tilde{P}+\tilde{Q}))^{\sigma}}{h(\tilde{P}+\tilde{Q})}&=&\dfrac{h(S+(\tilde{P}+\tilde{Q})^{\sigma}-(\tilde{P}+\tilde{Q}))}{h(S)}= \dfrac{h(S+(\tilde{P}+\tilde{Q})^{\sigma}-(\tilde{P}+\tilde{Q}))}{h(S+(\tilde{P})^{\sigma}-\tilde{P})} \dfrac{h(S+\tilde{P}^{\sigma}-\tilde{P})}{h(S)}\\ \\
		&=& 
		\dfrac{(h(\tilde{Q}))^{\sigma}}{h(\tilde{Q})} \dfrac{(h(\tilde{P}))^{\sigma}}{h(\tilde{P})}
	\end{eqnarray*}
	 This gives 
	\begin{equation*}
	\dfrac{h(\tilde{P}+\tilde{Q})}{h(\tilde{P})h(\tilde{Q})}=\dfrac{(h(\tilde{P}+\tilde{Q}))^{\sigma}}{(h(\tilde{P}))^{\sigma}(h(\tilde{Q}))^{\sigma}}=\left( \dfrac{h(\tilde{P}+\tilde{Q})}{h(\tilde{P})h(\tilde{Q})}\right)^{\sigma}
	\end{equation*}
	for every $ \sigma\in \Gal(\bar{K}/K) $. Now we have
	
	\begin{equation*}
	\dfrac{h(\tilde{P}+\tilde{Q})}{h(\tilde{P})h(\tilde{Q})}\in K, \textrm{ i.e.,  }\ \ h^{2}(\tilde{P}+\tilde{Q})\equiv h^{2}(\tilde{P})h^{2}(\tilde{Q}) \mod K^2.
	\end{equation*}
	Thus,
	\begin{equation*}
	g(P+Q)=g\circ[2](\tilde{P}+\tilde{Q})=h^{2}(\tilde{P}+\tilde{Q})\equiv h^{2}(\tilde{P})h^{2}(\tilde{Q})=g(P)g(Q) \mod K^2.
	\end{equation*}		
\end{Proof}

	\begin{Theorem}\label{thm1} Let $ C $ be a smooth genus 1 curve over $ \Q $ defined by an equation of the form
\begin{eqnarray*}
y^{2}=(a_{1}x+b_{1})(a_{2}x+b_{2})(a_{3}x+b_{3})(a_{4}x+b_{4}),\qquad \textrm{ where }a_i\in\Q^{\times},b_i\in \Q.
\end{eqnarray*}
Let $(x_0,y_0)\in C(\Q)$ be such that $x_0\ne-b_i/a_i$, $i=1,2,3,4$.
 We set $ \phi: C \longrightarrow E:=J(C) $ to be a $\Q$-birational isomorphism with $ \phi((x_0,y_0))=O_{E} $. For $Q\in C(\Q)$, one has
 
  $ Q\in 2C(\Q)$ if and only if
	$
	f_i(x_0)\,f_j(x_0)\,f_{i}(x(Q))f_{j}(x(Q))\in \Q^2
	$
 for all $ i,j \in\{1,2,3,4\} $ where $ f_{i}(x)=a_{i}x+b_{i} $.\\
\end{Theorem}
\begin{Proof}
The statement that $Q\in 2C(\Q)$ implies that $	f_i(x_0)\,f_j(x_0)\,f_{i}(x(Q))f_{j}(x(Q))\in \Q^2$ is a direct consequence of Proposition \ref{Prop1}.

So we assume that $ Q\in C(\Q) $ is such that $f_i(x_0)\,f_j(x_0)\,f_{i}(x(Q))f_{j}(x(Q))\in \Q^2$. Let $P\in  E(\Q) $ be such that $\phi^{-1}(P)=Q$. It suffices to show that if $ g_{ij}(P)\equiv 1 \mod  \Q^2 $ for all $ i,j $, then $ P\in 2E(\Q) $. 

Set $\textstyle [2]^{-1}P=\{R_{i}:1\le i\le 4 \}\subset E$. Let fix $ R\in [2]^{-1}P $. We also set $ R_{C,i}= \phi^{-1} (R_{i}) \in C $, $ i=1,2,3,4 $.  We recall that $R_i=R+T_i$ for some $T_i\in E[2]$. A simple calculation of the Jacobian $E$ shows that $E(\Q)[2]\isom\Z/2\Z\times\Z/2\Z$, hence all $2$-torsion points of $E$ are rational.  It follows that $[\Q(x(R_i),y(R_i)):\Q]$ is fixed for all $i=1,2,3,4$. 

Recall that $g_{ij}\in K(E)$, $ 1\le i, j\le 4$, is defined as follows 	
	\[g_{ij}(P)=
      f_i(x_0)f_j(x_0)f_i(x(\phi^{-1}(P)))f_j(x(\phi^{-1}(P)))
\] 
	where $x(\phi^{-1}(P))$ is the $x$-coordinate of $\phi^{-1}(P)\in C(K)$.  
By Proposition \ref{Prop2}, we have
\begin{equation*}
g_{ij}(P)=g_{ij}(2R_m)\equiv g_{ij}(R_m)^{2}\mod\Q^{2},\qquad m=1,2,3,4.
\end{equation*}
However, one knows that $g_{ij}(P)\in \Q^2$ by assumption. It follows that 
 $ g_{ij}(R_m)\in \Q $. Since $ g_{ij}(R_m)=f_i(x_0)f_j(x_0)f_{i}(x(R_{C,m}))f_{j}(x(R_{C,m})) $, it follows that $[\Q(x(R_{C,m})):\Q]\le2$. Writing 
$ x(R_{C,m})= A+B\sqrt{D} $ for some $A,B,D\in \Q$, one sees that 
\begin{equation*}
\dfrac{g_{12}(R_m)}{f_{1}(x_0)f_{2}(x_0)}=(a_1(A+B\sqrt{D})+b_{1})(a_{2}(A+B\sqrt{D})+b_{2}) \in \Q^*.
\end{equation*}
Therefore, either $ B=0 $ or $ a_{1}b_{2}+2a_{1}a_{2}A+a_{2}b_{1}=0 $. If $ B\neq 0 $, then $\textstyle A=\frac{-a_{1}b_{2}-a_{2}b_{1}}{2a_{1}a_{2}} $. In a similar fashion, since $\textstyle	\frac{g_{13}(R_m)}{f_{1}(x_0)f_{3}(x_0) } \in \Q$, one has 
$\textstyle A=\frac{-a_{1}b_{3}-a_{3}b_{1}}{2a_{1}a_{3}} $. One concludes that $ \frac{b_{2}}{a_{2}}=\frac{b_{3}}{a_{3}} $, which contradicts the fact that the points $(-b_i/a_i,0)\in C$ must be distinct. It follows that $ B=0 $,  i.e., $ x(R_{C,m})\in \Q$. Since $y(R_{C,m})^{2}=f_{1}(x(R_{C,m}))f_{2}(x(R_{C,m}))f_{3}(x(R_{C,m}))f_{4}(x(R_{C,m})) $, the latter implies that $ y(R_{C,m})\in \Q $ or $ y(R_{C,m})= K\sqrt{D} $ for some $ K\in \Q $ and $D\in \Q^*/\Q^{*2}$. In particular, $\Q(R_{C,m})= \Q(\sqrt{D})$. Since $\phi$ is a $\Q$-birational isomorphism, it follows that $\Q(R_m)= \Q(\sqrt{D})$. Moreover, from the observation above, all $R_m$ are $\Q$-rational, or all are defined over $\Q(\sqrt{D})\setminus\Q$.

One knows that since $P\in E(\Q)$, it follows that
\begin{equation*}
P=P^{\sigma}=(2R)^{\sigma}=2R^{\sigma}
\end{equation*}
for all $ \sigma \in \Gal(\overline \Q/\Q) $.  In addition, since $\phi^{-1}(S^{\sigma})=(\phi^{-1}(S))^{\sigma}$ for all $ \sigma \in \Gal(\overline \Q/\Q) $ and $S\in E$.
Therefore, we may assume without loss of generality that 
\begin{equation*}
R_{C,1}=(u_1, u_2\sqrt{D}), \ \  R_{C,2}=(u_1, -u_2\sqrt{D}), \ \  R_{C,3}=(v_1, v_2\sqrt{D}) \  \ R_{C,4}=(v_1, -v_2\sqrt{D}) ,
\end{equation*}
where $ u_1,v_1,u_2,v_2\in \Q$. 

We observe that $g_{ij}(R_{1})g_{ij}(R_{2})\in\Q^2$, we also see that $g_{ij}(R_{1})g_{ij}(R_{2})\equiv g_{ij}(P+T_1)$ for some $T_1\in E[2]$ for $i,j\in\{1,2,3,4\}$. In view of Proposition \ref{Prop2}, we obtain that $g_{ij}(T_1)\in\Q^2$ for $i,j\in\{1,2,3,4\}$. Repeating the argument above for $R_3$ and $R_4$, we get that $g_{ij}(T_2)$ for some $T_2\in E[2]$, $T_2\ne T_1$.  Noticing that $T_2\pm T_1\in E[2]$, it follows that $g_{ij}(T)\in\Q^2$ for all $i,j=1,2,3,4$, and all $T\in E[2]$. 

Using the fact that $g_{ij}(T)\in\Q^2$ for all $T\in E[2]$,  we may replace the point $P$ in the argument above with a point $T\in E[2]$.  In particular, as seen above,  this leads to the following:
Given $T_i\in E[2]$, $i=1,2,3,4$, there exists $T_i^j\in E$ such that $2T_i^j=T_i$, $j=1,2,3,4$, where $T_i^1,T_i^2,T_i^3,T_i^4$ are all in $E(\Q)$ or all in $E(\Q(\sqrt{d_i}))\setminus E(\Q)$ for some square free integer $d_i\ne0$.   Now since $2(T_i^j\pm T_s^t)\in E[2]$, then this implies that all $T_i^j\in E(K)$, for all $1\le i,j\le 4$, where $K$ is either $\Q(\sqrt{D})$ or $\Q(\sqrt{D},\sqrt{D'})$ where $D$ and $D'$ are square free integers. The fact that $T_{i}^j$ cannot be all in $E(\Q)$ is due to Mazur's classification of rational torsion points on elliptic curves defined over $\Q$, see \cite{Mazur1} or \cite[VIII.7, Theorem 7.5]{Silverman1}.

Now we rule out the possibility that $T_i^j$ are lying in $E(K)$ where $K$ is either $\Q(\sqrt{D})$ or $\Q(\sqrt{D},\sqrt{D'})$, hence $R_{C,i}$ should have all lived in $C(\Q)$ for $i=1,2,3,4$. If $K=\Q(\sqrt{D})$, then this means that the torsion part of $E(K)$ contains complete $2$-torsion where for each $2$-torsion point $T_i$ there are $4$ distinct torsion point $T_{i}^j$ such that $4T_i^j=T_i$. A complete classification of possible torsion points of elliptic curves over quadratic fields was established in \cite{Kam,Kenku} after a series of papers. In particular, if an elliptic curve is defined over $E(\Q(\sqrt{D}))$ with a non-cyclic torsion group containing $\Z/2\Z\times\Z/2\Z$, then it should be one of the following groups 
\[\Z/2\Z\times\Z/2m\Z,1\le m\le 6,\qquad\textrm{or } \Z/4\Z\times\Z/4\Z\]
where it can be seen easily that it is impossible for the $2$-torsion points to satisfy the aforementioned property. 

Now we rule out the possibility that $K=\Q(\sqrt{D},\sqrt{D'})$.  In \cite{Chou}, a complete classification for torsion points of elliptic curves defined over Galois quartic fields is given. If an elliptic curve is defined over $E(K)$ with a non-cyclic torsion group containing $\Z/2\Z\times\Z/2\Z$, then it should be one of the following groups  \cite[Theorem 1.4]{Chou}
\[\Z/2\Z\times\Z/2m\Z,1\le m\le 8,\qquad\Z/4\Z\times\Z/4n\Z, n=1,2, \qquad\textrm{or }\Z/6\Z\times\Z/6\Z.\]
Again one may check that for neither of these groups all two torsion points are divisible by $4$. Hence when $g_{ij}(P)\in\Q^2$, for all $ i,j=1,2,3,4$, then $R_{C,i}\in C(\Q)$, in particular $Q=\phi^{-1}(P)\in 2C(\Q)$.
\end{Proof}

 \begin{Remark}
	\label{Rem1}
	In Theorem \ref{thm1}, if $(x_0,y_0)$ is chosen to be $\infty_+$, then $g_{ij}$ becomes  
	$$g_{ij}(P)=	a_i a_j f_{i}(x(\phi^{-1}(P)))f_{j}(x(\phi^{-1}(P))).$$ In this case we assume $a_1a_2a_3a_4\in (K^{\times})^2$ to make sure that the points $\infty_{\pm}$ are in $C(\Q)$.
\end{Remark}

\section{Examples}
\label{sec:examples}
 
 Let $C$ be defined over a perfect field $K$ by the following quartic model
 	\begin{eqnarray*}
 		v^{2}=au^{4}+bu^{3}+cu^{2}+du+q^{2}
 	\end{eqnarray*}
 	where $a,b,c,d \in K$ and $q\in K^{*}$. The birational isomorphism $\phi_1 :C\to E:=J(C)$ is defined by 
 	\begin{eqnarray}\label{birational}
 		x= (2q(v+q)+du)/u^{2} \quad \textrm{and} \quad y=(4q^{2}(v+q)+2q(du+cu^{2})-d^{2}u^{2}/2q)/u^{3}
 	\end{eqnarray}
 	where $E$ is described by \[y^2+a_1xy+a_3y=x^3+a_2x^2+a_4x+a_6,\]
 	and 
 	\begin{eqnarray*}
 		a_{1}=d/q,\ \ \ a_{2}=c-d^{2}/4q^{2}, \ \ \ a_{3}=2qb, \ \ \ a_{4}=-4q^{2}a, \ \ \ a_{6}= a_{2}a_{4}, 
 	\end{eqnarray*}
 	see \cite[Chapter 1, Proposition 1.2.1]{Iann}. The inverse map is given by
 	\begin{eqnarray*}
 		u=(2q(x+c)-d^{2}/2q)/y, \ \ \ \ \ \ v= -q+u(ux-d)/2q.
 	\end{eqnarray*}
 
In view of Theorem \ref{thm1}, the point $(x_0,y_0)$ is $(0,q) $ as $\phi(0,q)=O_E$. 

\begin{Example}
	   Let $C$ : \  $v^{2}=(u+1)(4u+1)(5u+1)(14u+1)$. We consider the map $\phi_1:C\to E$ defined above,  where
	\begin{eqnarray*}
		E : \ \ y^{2}+24xy+852y=x^{3}+25x^{2}-1120x-28000.
	\end{eqnarray*}
	
	We start with the point $T=(1,30)\in C(\Q)$. Then
	\begin{eqnarray*}
		\phi_1(T)=(86,222), \quad\textrm{and }  2\phi_1(T)=(-14335/784, -289575/21952)\in E(\Q).
	\end{eqnarray*}
	Now $Q=\phi_1^{-1}(2\phi_1(T))=(-56/55, 1053/605)\in 2C(\Q)$. Setting $f_1(u)=u+1$, $f_2(u)=4u+1$, $f_3(u)=5u+1$ and $f_4(u)=14u+1$,  we see that
	\begin{eqnarray*}
		f_{1}\left(\frac{-56}{55}\right)=\frac{-1}{55}, \ \ f_{2}\left(\frac{-56}{55}\right)= \frac{-169}{55}, \ \ f_{3}\left(\frac{-56}{55}\right)=\frac{-225}{55}, \ \ f_{4}\left(\frac{-56}{55}\right)=\frac{-729}{55}.
	\end{eqnarray*}
Since $f_i(0)=1$ for $i=1,2,3,4$, it follows that \[g_{ij}(2\phi_1(T))=f_i\left(-56/55\right)f_j\left(-56/55\right)\in \Q^2\quad\textrm{for all }  i,j=1,2,3,4.\] 
\end{Example}

\begin{Example}
	Let $C$ : \ $v^2=(u+2)(u+4)(u+8)(u+9)$. Then the transformation $\phi_1$ gives the elliptic curve
	\begin{eqnarray*}
		E : \ \ y^{2}+\frac{71}{3}xy+1104y=x^{3}+\frac{1511}{36}x^{2}-2304x-96704.
	\end{eqnarray*}
	Considering the point $T=(-5,6)\in C(\Q)$,  we obtain
	\begin{eqnarray*}
		\phi_1(T)=(-56,404/3), \quad\textrm{and }  2\phi_1(T)=(-199/4, 245/6)\in E(\Q).
	\end{eqnarray*}
	Now the point $\phi_1^{-1}(2\phi_1(T))=(-64/7, -120/49)$ lies in $2C(\Q)$. We see that
	
	\begin{eqnarray*}
		f_{1}\left(\frac{-64}{7}\right)=\frac{-50}{7}, \ \ f_{2}\left(\frac{-64}{7}\right)= \frac{-36}{7}, \ \ f_{3}\left(\frac{-64}{7}\right)=\frac{-8}{7}, \ \ f_{4}\left(\frac{-64}{7}\right)=\frac{-1}{7}
	\end{eqnarray*}
	where $f_1(u)=u+2$, $f_2(u)=u+4$, $f_3(u)=u+8$ and $f_4(u)=u+9$.
	Since $ f_{1}(0)=2, \ f_{2}(0)=4, \ f_{3}(0)=8, \ f_{4}(0)=9 $,  we see that $$g_{ij}(2\phi_1(T))=f_{i}(0)f_j(0)f_i(-64/7)f_j(-64/7)\in\Q^2\quad\textrm{ for all }i,j.$$ 
\end{Example}

 In the following example, we consider a birational map between a quartic model and its jacobian elliptic curve different from the map introduced in the previous two examples. The elliptic curve
  \begin{eqnarray*}
		E : w^{2}=v^{3}+Av+B,
	\end{eqnarray*}
	where $  A,B \in K $, is the Jacobian of the elliptic curve defined by the following quartic model \begin{eqnarray*}
		C: y^{2}=x^{4}-6ax^{2}-8bx+c
	\end{eqnarray*}
	where $ c=-4A-3a^{2} $ and $ B=b^{2}-a^{3}-Aa .$  We notice that the point $ P=(a,b)\in E(K) $ and that $\infty_+$ and $\infty_-$ map to $O_E$ and $P$ respectively, see \cite[\S 2]{Adams}. We define a birational isomorphism  $\phi_2:C\to E$ as follows
	\begin{eqnarray*}
		x=\frac{w+b}{v-a}, \ \ \ \ y=2v+a-\left(\frac{w+b}{v-a}\right)^{2},
	\end{eqnarray*} 
	whereas the inverse map of $\phi_2$ is given by
	\begin{eqnarray*}
		v=\frac{1}{2}(x^{2}+y-a), \ \ \ \ w=\frac{1}{2}(x^{3}+xy-3ax-2b).
	\end{eqnarray*}

\begin{Example}
 Let $ C: y^2=(x+2)(x+4)(x+8)(x+9) $. Applying the transformation $x\mapsto x-\frac{23}{4}$, $y\mapsto y$, we get the quartic curve
	\begin{eqnarray*}
	\tilde{C}: y^2= x^4-\frac{131}{8}x^2-\frac{33}{8}x+\frac{12285}{256}.
	\end{eqnarray*}
Setting $ a=\frac{131}{48} , \ b=\frac{33}{64}, \ \ A=-\frac{211}{12}, \ \ B=\frac{754}{27}$, the Jacobian elliptic curve is defined by
\begin{eqnarray*}
E: w^2=v^3-\frac{211}{12}v+\frac{754}{27}
\end{eqnarray*}

The point $Q=(-5,6)$ is in $ C(\Q) $ and $ 2Q=(-\frac{217}{24},\frac{715}{576}) $. We get
\begin{eqnarray*}
	f_{1}\left(-\frac{217}{24}\right)=-\frac{169}{24}, \ \ f_{2}\left(-\frac{217}{24}\right)=-\frac{121}{24}, \ \ f_{3}\left(-\frac{217}{24}\right)=-\frac{25}{24},\ \ f_{4}\left(-\frac{217}{24}\right)=-\frac{1}{24}\\
\end{eqnarray*}
where $f_1(x)=x+2$, $f_2(x)=x+4$, $f_3(x)=x+8$ and $f_4(x)=x+9$.
 According to Remark \ref{Rem1}, $a_i=1$ for all $i=1,2,3,4$. It follows that $$g_{ij}(2\phi_1(Q))=f_i(-217/24)f_j(-217/24)\in\Q^2\quad\textrm{ for all }i,j.$$ 
\end{Example}

\section{$4$-torsion points on quartic models}

In what follows we give a necessary and sufficient condition for an elliptic curve defined by a quartic model over $\Q$ to possess a $4$-torsion point. 

\begin{Theorem}
\label{thm:tor4}
	Let $C$ be a smooth genus one curve defined over $  \Q$ by  
	\begin{equation*}
	v^2=(k_{1}u+1)(k_{2}u+1)(k_{3}u+1)(k_{4}u+1),\qquad k_i\in\Q^{\times}.
	\end{equation*}
	Fix the map $\phi_1:C\to E:=J(C)$ defined as in \S \ref{sec:examples}.
 Then $C$ has a $4$-torsion point defined over $\Q$ if and only if one of the following holds:
\begin{itemize} 
 \item[i)] $ (k_{1}-k_{3})(k_{2}-k_{4}) $ and $ (k_{1}-k_{2})(k_{3}-k_{4}) $ are both squares in $\Q$; or
\item[ii)] $(k_{2}-k_{3})(k_{1}-k_{4}) $ and $ (k_{1}-k_{2})(k_{4}-k_{3}) $ are both squares in $\Q$; or
\item[iii)] $(k_{3}-k_{2})(k_{1}-k_{4}) $ and $ (k_{1}-k_{3})(k_{4}-k_{2}) $ are both squares in $\Q$.
\end{itemize}
\end{Theorem}
\begin{Proof}
	The curve $ J(C) $ is defined by the Weierstrass equation
	\begin{equation*}
E:= y^2+a_{1}xy+a_{3}y=x^3+a_{2}x^2+a_{4}x+a_{6}
	\end{equation*}
where $ a_{1}=k_{1} + k_{2} + k_{3} + k_{4}, \  a_{2}=k_{1} k_{2} + k_{1} k_{3} + k_{2} k_{3} + k_{1} k_{4} + k_{2} k_{4} + k_{3} k_{4} - 
	\frac{1}{4} (k_{1} + k_{2} + k_{3} + k4)^2, \ a_{3}= 2 (k_{1} k_{2} k_{3} + k_{2} k_{3} k_{4} + k_{1} (k_{2} + k_{3}) k_{4}), \ a_{4}= -4 k_{1} k_{2} k_{3} k_{4} , \ a_{6}= a_{2}a_{4}.$	
	Using MAGMA, the $2$-torsion points on $ E $ are
	
	\begin{eqnarray*}
	P_{1}&=& O_{E}, \\ 
	P_{2}&=& (-k_{1} k_{4} - k_{2} k_{3}, \ \frac{1}{2} k_{1}^2 k_{4} - \frac{1}{2} k_{1} k_{2} k_{3} - \frac{1}{2} k_{1} k_{2} k_{4} - \frac{1}{2} k_{1} k_{3} k_{4} + \frac{1}{2} k_{1} k_{4}^2 + \frac{1}{2} k_{2}^2 k_{3} + \frac{1}{2} k_{2} k_{3}^2 - \frac{1}{2} k_{2} k_{3} k_{4}), \\ 
	P_{3}&=& (-k_{1} k_{3} - k_{2} k_{4}, \  \frac{1}{2} k_{1}^2 k_{3} - \frac{1}{2} k_{1} k_{2} k_{3} - \frac{1}{2} k_{1} k_{2} k_{4} + \frac{1}{2} k_{1} k_{3}^2 - \frac{1}{2} k_{1} k_{3} k_{4} + \frac{1}{2} k_{2}^2 k_{4} - \frac{1}{2} k_{2} k_{3} k_{4} + \frac{1}{2} k_{2} k_{4}^2), \\ 
	P_{4}&=& ( -k_{1} k_{2} - k_{3} k_{4}, \frac{1}{2} k_{1}^2 k_{2} + \frac{1}{2} k_{1} k_{2}^2 - \frac{1}{2} k_{1} k_{2} k_{3} - \frac{1}{2} k_{1} k_{2} k_{4} - \frac{1}{2} k_{1} k_{3} k_{4} - \frac{1}{2} k_{2} k_{3} k_{4} + \frac{1}{2} k_{3}^2 k_{4} + \frac{1}{2} k_{3} k_{4}^2).
	\end{eqnarray*}
Then 
\begin{eqnarray*}
 \phi_1^{-1}(P_2)&=&((k_{1} - k_{2} - k_{3} + k_{4})/(k_{2} k_{3} - k_{1} k_{4}), -(((k_{1} - k_{2}) (k_{1} - k_{3}) (k_{2} - k_{4}) (k_{3} - k_{4}))/(k_{2} k_{3} - k_{1} k_{4})^2)),\\
 \phi_1^{-1}(P_{3})&=& ( (-k_{1} + k_{2} - k_{3} + k_{4})/(k_{1} k_{3} - k_{2} k_{4})
, ((k_{1} - k_{2}) (k_{2} - k_{3}) (k_{1} - k_{4}) (k_{3} - k_{4}))/(k_{1} k_{3} - k_{2} k_{4})^2),\\
 \phi_1^{-1}(P_{4})&=& ((-k_{1} - k_{2} + k_{3} + k_{4})/(k_{1} k_{2} - k_{3} k_{4}), ((k_{1} - k_{3}) (-k_{2} + k_{3}) (k_{1} - k_{4}) (k_{2} - k_{4}))/(k_{1} k_{2} - k_{3} k_{4})^2  ).
\end{eqnarray*}
In view of Theorem \ref{thm1}, the point $ \phi_1^{-1}(P_2)=(u_2,v_2)\in 2C(\Q)$ if and only if 
$ (k_{i}u_{2}+1)(k_{j}u_{2}+1) \in \Q^{2}$, where $ i,j\in\{1,2,3,4\} $.  Now direct substitution yields that the latter conditions are equivalent to $(k_{1}-k_{3})(k_{2}-k_{4}) $ and $ (k_{1}-k_{2})(k_{3}-k_{4}) $ are squares in $\Q$. The other two conditions follow by considering the points  $\phi_1^{-1}(P_{3}) $ and  $\phi_1^{-1}(P_{4})$.  
\end{Proof}

As an application of Theorem \ref{thm:tor4}, we construct the following example of a $2$-parameter family of elliptic curves over $\Q$ described by a quartic equation for which none of the nonsingular fibers has a nontrivial rational torsion point of order $4$.

\begin{Example}
 Let $ C_{s,t} $ be a smooth genus 1 curve over $\Q(s,t) $ defined by 
 \begin{equation*}
v^2 =\left(\frac{t^2+s^2+ts}{t+s}u+1\right)\left(\frac{-ts}{t+s}u+1\right)(tu+1)(su+1).
 \end{equation*}

	Setting $ k_{1}=\frac{t^2+s^2+ts}{t+s}, \ k_{2}= \frac{-ts}{t+s} $, one may use Theorem \ref{thm:tor4} to investigate the existence of rational numbers $s,t$ such that $C_{s,t}$ has a torsion point of order $4$. 
Condition (iii) of Theorem \ref{thm:tor4} is the fact that the two expressions $ t(2s+t) $ and $ s(2t+s) $ are squares in $\Q$. The latter is equivalent to the existence of a rational point on the following intersection of two quadric surfaces in $\PP^3$:\[u^2=t(2s+1),\qquad v^2=s(2t+s).\]
 In fact, the latter intersection is an elliptic curve that can be described by the Weierstrass equation $y^2 = x^3 - 4x^2 + 16x$ and whose Mordell-Weil group is isomorphic to $\Z/4\Z$ corresponding to the points $(s:t:u:v)=(0:1:\pm1:0), (1:0:0:\pm1)$.  
 
 Conditions (i) and (ii) of Theorem \ref{thm:tor4} are equivalent to the existence of a rational point on the intersection of the quadric surfaces 
\begin{eqnarray*}
u^2&=& -s(s+2t) ,\qquad  v^2=-s^2+t^2,\textrm{ and} \\
u^2&=&-t(t+2s), \qquad v^2=-t^2+s^2,
\end{eqnarray*} respectively.  Both intersections are isomorphic to the elliptic curve described by $y^2=x^3-x^2+x$ whose Mordell-Weil group is isomorphic to $\Z/4\Z$ corresponding to the points $(s:t:u:v)=(0:1:0:\pm1), (1:-1:\pm1:0)$. 

 It follows that the curve $C_{s,t}$ does not have a torsion point of order $4$ over $\Q$ for any choice of the rational pair $s,t$ with $s,t\ne 0,s\ne \pm t$.
\end{Example}

\section{Diophantine $D(q)$-quintuples}

\begin{Definition}
Let $q$ be a nonzero rational number. A set of $m$ nonzero rational numbers $\{a_1, a_2, \cdots , a_m\}$ is called a {\em rational $D(q)$-$m$-tuple} if $a_i  a_j + q$ is a perfect rational square for all $1 \le i < j \le m$. If $q=1$, then this set is called a {\em rational Diophantine $m$-tuple}.
\end{Definition}

There are only finitely many ways of extending a rational
$D(q)$-quadruple to a rational $D(q)$-quintuple, see \cite{Herrmann}. In \cite{ondiophantine}, an explicit expression for the element extending a rational $D(q)$-quadruple to a raional $D(q)$-quintuple was provided if $q$ is a rational square.

The definition above can be extended over the ring of polynomial with integer coefficients as follows.  

\begin{Definition}
Let $q\in\Q[x]$ be a nonzero polynomial. Let $\{a_1, a_2,\cdots , a_m\}$ be a set of $m$ nonzero polynomials with rational coefficients. We assume that there does not exist a polynomial $p\in \Q[x]$ such that $a_1/p,\cdots, a_m/p$ and $q/p^2$ are rational numbers. The set $\{a_1, a_2,\cdots , a_m\}$ is called a {\em polynomial $D(q)$-$m$-tuple} if  $a_ia_j+q=b_{ij}^2$, for all $1\le i < j\le m$, where $b_{ij}\in\Q[x]$.
\end{Definition}
 The assumption that there is no polynomial $p$ such that $a_1/p,\cdots , a_m/p$ and $q/p^2$ are rational numbers implies that if $q$ is constant then not all elements $a_1,\cdots , a_m$ of a polynomial $D(q)$-$m$-tuple are allowed to be constant. When $q$ is a linear polynomial, the latter condition is trivially always satisfied.

In what follows we will be interested in polynomial $D(q)$-$m$-tuples whose elements are linear polynomials, and $q$ is also a linear polynomial. We define 
$$L_{1}=\sup\{|S|: S \textrm{ is a polynomial $D(ax+b)$-tuple consisting of linear polynomials for some }a\ne 0 \textrm{ and }b\}.$$ It was shown that $L_{1}=4$, see \cite[Theorem 1]{Dujella}.  Therefore, the set $\{x,16x+8,25x+14,36x+20\} $ which is a polynomial $D(16x+9)$-quadruple cannot be extended to a polynomial $D(16x+9)$-quintuple using a linear polynomial.
In this section, we show that for infinitely many rational values of $x$, the latter rational $D(16x+9)$-quadruple can be extended to a rational $D(16x+9)$-quintuple. The main tool is the following straightforward corollary of Theorem \ref{thm1}.

\begin{Corollary}
\label{cor1}
Let $ C $ be a smooth genus 1 curve over $ \Q $ defined by an equation of the form
\begin{eqnarray*}
y^{2}=(a_{1}x+b_{1})(a_{2}x+b_{2})(a_{3}x+b_{3})(a_{4}x+b_{4}),\qquad \textrm{ where }a_i\in\Q^{\times},b_i\in \Q.
\end{eqnarray*}
 We set $ \phi: C \longrightarrow E:=J(C) $ to be a $\Q$-birational isomorphism with $ \phi((x_0,y_0))=O_{E} $ for some $(x_0,y_0)\in C(\Q)$, $x_0\ne -b_i/a_i$, $i=1,2,3,4$.  Assume, moreover, that $f_i(x_0)f_j(x_0)\in \Q^2$, for all $i,j=1,2,3,4$, where $f_i(x)=a_i x+b_i$.  Let $Q\in C(\Q)$.
 
 Then $ Q\in 2C(\Q)$ if and only if there exists $\delta_Q\in \Q$ such that
	$$ a_{i}x(Q)+b_{i}=\delta_Q \cdot z_i^2 \textrm{ for some }z_i\in\Q $$
 for all $ i \in\{1,2,3,4\} $.
\end{Corollary}

Corollary \ref{cor1} implies the following result on extending a rational $D(q)$-quadruple to a $D(q)$-quintuple. 

\begin{Corollary}
\label{cor2}
Let $q$ be a nonzero rational number. Let $S=\{a_1,a_2,a_3,a_4\}$ be a $ D(q) $-quadruple.  
Consider the smooth genus one curve $$C_S:y^2=(a_1x+q)(a_2x+q)(a_3x+q)(a_4x+q).$$
We fix the birational isomorphism $ \phi_1: C_S \longrightarrow J(C_S) $ defined as in \S \ref{sec:examples}. Then $S$ can be extended to a rational $D(q)$-quintuple if and only if there is a point $Q\in 2C_S(\Q)$ with $\delta_Q\in\Q^2$.  
\end{Corollary} 

In the following theorems, we extend some of the known polynomial $D(q)$-quadruples consisting of linear polynomials to $D(q)$-quintuples for infinitely many rational values of $x$. 

\begin{Theorem}\label{thm2}
	The polynomial $ D(16t+9) $-quadruple $\{t, 16t+8,2 25t+14, 36t+20 \}$ can be extended to a rational $ D(16t+9) $-quintuple for infinitely many $ t\in \mathbb{Q} $. 
\end{Theorem}
\begin{Proof}
Let 
$ k_{1}= t $,
$ k_{2}= 16t+8 $,
$ k_{3}= 25t+14 $,
$ k_{4}= 36t+20 $, 
$ q= 16t+9 $. Consider the smooth genus one curve $ C $ defined by
\begin{equation}\label{quarticc}
 v^2= (k_1 u + q) (k_2 u + q) (k_3 u + q) (k_4 u + q)
 \end{equation}
  over $\Q(t)$. By using the birational isomorphism $ \phi_1:C\to E $, defined in \S \ref{sec:examples}, we have the elliptic curve
\begin{equation*}
E: y^{2}+a_{1}xy+a_{3}y=x^{3}+a_{2}x^{2}+a_{4}x+a_{6}
\end{equation*} 
with the coefficients 
\begin{eqnarray*}
 a_{1} &=& 6(7+13t)(9+16t),\\
   a_{2}&=&(9+16 t)^{2}(111+8t(55+54t)),\\
  a_{3} &=& 8(7+13t) (9 + 16 t)^3 (80 + t (318 + 313 t)), \\
   a_{4} &=& -128 t (1 + 2 t) (5 + 9 t) (9 + 16 t)^4 (14 + 25 t), \\
    a_{6} &=& -128 t (1 + 2 t) (5 + 9 t) (9 + 16 t)^6 (14 + 25 t) (111 + 
8 t (55 + 54 t))
\end{eqnarray*}
where $(x_0,y_0)=(0,q^2)$ is such that $\phi_1(x_0,y_0)=O_E$.

We take the point $ P=(0, -q^{2})\in C(\Q(t))$. Then
 \begin{equation*}
 S:=\phi_1(P)= (-(9 + 16 t)^2 (111 + 8 t (55 + 54 t)), \ 2 (1 + 2 t) (7 + 13 t) (9 + 16 t)^3 (13 + 22 t)).
 \end{equation*} 
 
  Using MAGMA, \cite{Magma}, $ Q=2S = (x_{1}, y_{1}) $ is given by
\begin{equation*}
x_{1}= \frac{u(t)}{(4 (1 + 2 t)^2 (7 + 13 t)^2 (13 + 22 t)^2)},\quad y_{1}= -\frac{v(t)}{8 (1 + 2 t)^3 (7 + 13 t)^3 (13 + 22 t)^3}.
\end{equation*}
where
\begin{eqnarray*}
\lefteqn {u(t)= 4965468561 + 87791232672 t + 698248164432 t^2 + 3289862320448 t^3 + 
	10168707377552 t^4 + 21544947073664 t^5 + {} }
\nonumber \\
& & {}+ 31689009677248 t^6 + 
	31949101618688 t^7 + 21130944883712 t^8 + 8278920101888 t^9 + 
	1459071221760 t^{10}
	\end{eqnarray*}

and

{\footnotesize$$v(t)=(9 + 16 t)^4 (85 + 302 t + 268 t^2) (111 + 398 t + 356 t^2) (71 + 
	8 t (31 + 27 t)) (97 + 8 t (43 + 38 t)) (1369 + 
	32 t (232 + t (419 + 252 t))).$$}

Then $ \phi_1^{-1}(Q)=(u_{1}, v_{1}) $, where

\begin{equation*}
u_{1}=-\frac{4 (1 + 2 t) (7 + 13 t) (13 + 22 t)}{1369 + 32 t (232 + t (419 + 252 t))},
\end{equation*}\\

\begin{equation*}
v_{1}= \frac{(85 + 302 t + 268 t^2) (111 + 398 t + 356 t^2) (71 + 
	8 t (31 + 27 t)) (97 + 8 t (43 + 38 t))}{(1369 + 32 t (232 + t (419 + 252 t)))^2}.
\end{equation*}\\

Therefore, we obtain

 \begin{equation*}
 k_{1}u_{1}+q=\frac{(111 + 398 t + 356 t^2)^2}{1369 + 32 t (232 + t (419 + 252 t))},
 \end{equation*}\\
 \begin{equation*}
 k_{2}u_{1}+q= \frac{(97 + 8 t (43 + 38 t))^2}{1369 + 32 t (232 + t (419 + 252 t))},
 \end{equation*}\\
\begin{equation*}
k_{3}u_{1}+q= \frac{(85 + 302 t + 268 t^2)^2}{1369 + 32 t (232 + t (419 + 252 t))},
\end{equation*}\\
\begin{equation*}
k_{4}u_{1}+q= \frac{(71 + 8 t (31 + 27 t))^2}{1369 + 32 t (232 + t (419 + 252 t))}.
\end{equation*}\\
One has $\delta_Q=\frac{1}{1369 + 32 t (232 + t (419 + 252 t))}.$ In view of Corollary \ref{cor2}, if $\delta_Q\in \Q^2$, then $k_iu_1+q\in\Q^2$ for all $i=1,2,3,4$. 

The elliptic curve $ r^2=1369 + 32 t (232 + t (419 + 252 t)) $ has Mordell-Weil rank $2$ over $ \mathbb{Q} $, \cite{Magma}. It follows that there are infinitely many  $ t\in \mathbb{Q} $  such that $\delta_Q\in\Q^2$, and hence the set $$\left\{t, 16t+8,2 25t+14, 36t+20, -\frac{4 (1 + 2 t) (7 + 13 t) (13 + 22 t)}{1369 + 32 t (232 + t (419 + 252 t))} \right\}$$ is a rational $ D(16t+9) $-quintuple.

\end{Proof}

By choosing $t$ to be the $t$-coordinate of a rational point on the elliptic curve $ r^2=1369 + 32 t (232 + t (419 + 252 t)) $,  the set $$\left\{t, 16t+8,2 25t+14, 36t+20, -\frac{4 (1 + 2 t) (7 + 13 t) (13 + 22 t)}{1369 + 32 t (232 + t (419 + 252 t))} \right\}$$ is a rational $D(q)$-quintuple produced by extending the polynomial $D(q)$-quadruple in Theorem \ref{thm2} when evaluated at $t$.  In the following table, we give examples of such $D(q)$-quintuples. 
\vskip5pt
\begin{center}
\begin{tabular}{|p{3.7cm}|c|c|}
\hline
$t$& $q$ & $D(q)$-quintuple\\
	\hline
 $ t= \frac{672}{8064} $ & $ q=\frac{31}{3} $ & $ \{\frac{1}{2}, \frac{28}{3}, \frac{193}{12}, 23, -\frac{60431}{225228}\} $ \\
	\hline
	 $ t= -\frac{3264}{8064} $ & $ q=\frac{53}{21} $ & $ \{-\frac{17}{42}, \frac{32}{21}, \frac{163}{42}, \frac{38}{7}, -\frac{50224}{240429}\} $ \\
	\hline
	$ t=\frac{-3600}{8084} $ & $ q=\frac{13}{7} $ & $ \{-\frac{25}{56}, \frac{6}{7}, \frac{159}{56}, \frac{55}{14}, -\frac{17889}{103544}\} $ \\
	\hline
	$  t=\frac{-4192}{8084} $ & $ q=\frac{43}{63} $ &  $ \{-\frac{131}{252}, -\frac{20}{63}, \frac{253}{252}, \frac{9}{7}, \frac{60085}{183708}\} $ \\
	\hline
	 $ t= -\frac{4572}{8064} $ & $ q=-\frac{1}{14} $ & $ \{-\frac{127}{224}, -\frac{15}{14}, -\frac{39}{224}, -\frac{23}{56}, -\frac{73455}{123704}\} $ \\
	 \hline
	  $ t= -\frac{4615}{8064} $ & $ q=-\frac{79}{504} $ & $ \{-\frac{4615}{8064}, -\frac{583}{504}, -\frac{2479}{8064}, -\frac{135}{224}, -\frac{3414104551}{6009297336}\} $\\
	 \hline

\end{tabular}
\end{center}

\begin{Theorem}\label{thm3} The following polynomial $ D(q) $-quadruples can be extended to a rational $ D(q) $-quintuple for infinitely many $ t\in \mathbb{Q} $. 
	\begin{itemize}
	\item[(i)] $ \{4t, 144t+8, 25t+1, 49t+3\} $, where $q=16t+1$,  can be extended using $ \frac{-4 (2 + 37 t) (3 + 58 t) (5 + 82 t)}{-1 + 32 t (13 + t (529 + 5148 t))} $, where $ t $ is the $ t $-coordinate of a rational point on the elliptic curve $ E: r^2= -1 + 32 t (13 + t (529 + 5148 t)) $.
	\item[(ii)] $ \{t, 9t+26, 4t+12, 16t+40\} $, where $q=16t+49$, can be extended using $ \frac{4 (1 + 2 t) (13 + 5 t) (27 + 10 t) (49 + 16 t)}{96721 + 16 t (6521 + 2342 t + 280 t^2)} $, where $ t $ is the $ t $-coordinate of a rational point on the elliptic curve $E: r^2= (96721 + 16 t (6521 + 2342 t + 280 t^2)) (16 t + 49) $.
	\item[(iii)] $ \{t, \frac{t}{4}-1, \frac{9t}{4}+5, 4t+8\} $, where $q=4t+9$, can be extended using $ \frac{(2 + t) (9 + 4 t) (8 + 5 t) (14 + 5 t)}{324 + 8 t (62 + t (31 + 5 t))} $, where $ t $ is the $ t $-coordinate of a rational points on the elliptic curve $E: r^2= (81 + 2 t (62 + t (31 + 5 t)))(9 + 4 t) $.
\end{itemize}
\end{Theorem}

\begin{Proof}
The proof is similar to to the proof of Theorem \ref{thm2}.
\end{Proof}
\begin{Remark}
The Mordell-Weil rank $r_E$ of $E(\Q)$ in Theorem \ref{thm3} is $r_E=3$ in (i) and (ii), whereas $r_E=2$ in (iii).
\end{Remark}

\begin{Example}
	Let $ t=\frac{180}{121} $ in the $ D(4t+1) $-quadruple given in Theorem \ref{thm3}.  We then obtain the following $ D(\frac{841}{121}) $-quintuple
	
	\begin{equation*}
	\left\{ \frac{180}{121}, -\frac{318}{121}, \frac{284}{121}, -\frac{248}{121}, \frac{2562308340}{2164017361}\right\}.
	\end{equation*}
	On the other hand, Theorem 1  in  \cite{ondiophantine} gives us that if $q, x_{1}, x_{2}, x_{3}, x_{4}$ are rational numbers such that $x_{i}x_{j} +q^2 = y_{ij}^{2}, \ \ y_{ij}\in \mathbb{Q}$, for all $1 \leq i < j \leq 4$ and $x_1x_2x_3x_4 \neq q^4$, then $ \{x_{1}, x_{2}, x_{3}, x_{4}, x_{5}\} $ is a rational $ D(q^2) $-quintuple, where
	 $x_{5} = A/B $ with 
	 
	 \begin{eqnarray*}
	 A&=& q^{3}\Big( \pm 2y_{12}y_{13}y_{14}y_{23}y_{24}y_{34} + qx_{1}x_{2}x_{3}x_{4}(x_{1} + x_{2} + x_{3} + x_{4})+
	 {}
	 \nonumber\\
	 && {}+2q^{3}(x_{1}x_{2}x_{3} + x_{1}x_{2}x_{4} + x_{1}x_{3}x_{4} + x_{2}x_{3}x_{4}) + q^{5}(x_{1} + x_{2} + x_{3} + x_{4}) \Big),\\
	 B &=& (x_{1}x_{2}x_{3}x_{4} - q^{4})^{2}.
	  \end{eqnarray*}
 Therefore, the rational $D(\frac{841}{121})$-quadruple $\{ 180/121, -318/121, 284/121, -248/121\}$ can be extended to a $D(\frac{841}{121})$-quintuple using either the rational number $ \frac{1255545720}{540051121} $ or $-\frac{143212695780}{74048750161}$. In particular, we obtain two almost $ D(\frac{841}{121}) $-sextuple, i.e., $x_ix_j+q^2$ is a rational square for all $1\le i<j\le 6$ except when $(i,j)=(5,6)$.   
 
 \end{Example}

\end{document}